\documentclass{amsart} 
\usepackage{amssymb}
\usepackage{mathrsfs}
\usepackage{cite}
\usepackage{graphicx}
\usepackage{caption}
\usepackage{subcaption}
\usepackage{float}
\usepackage{color}
\usepackage[foot]{amsaddr}

\newtheorem{thm}{Theorem}[section]

\theoremstyle{definition}
\newtheorem{definition}[thm]{Definition}

\theoremstyle{remark}
\newtheorem{remark}[thm]{Remark}

\numberwithin{equation}{section}

\DeclareMathOperator{\Real}{\mathbb{R}}  
\DeclareMathOperator{\ds}{d\it{s}}          
\DeclareMathOperator{\dt}{d\it{t}}                 
\DeclareMathOperator{\dx}{d\it{x}}                 
\DeclareMathOperator{\dy}{d\it{y}}                 

\DeclareMathOperator{\la}{\lambda}

\DeclareMathOperator{\ddt}{\frac{d}{d\textit{t}}} 

\begin{document}

\title[Regularity Criterion for 3D Boussinesq Equations]{Regularity Criterion for the Three-dimensional Boussinesq Equations}
\author{Karen Zaya}
\address{Department of Mathematics. University of Michigan.  2074 East Hall. 530 Church Street. Ann Arbor, MI 48109-1043}
\email{zaya@umich.edu}

\keywords{Boussinesq equations, regularity, dissipation wave number}
\subjclass[2010]{76D03, 35Q35}

\begin{abstract}
We prove that a solution $(u, \theta)$ to the three-dimensional Boussinesq equations does not blow-up at time T if $\| u_{\le Q}\|_{B^1_{\infty, \infty}}$ is integrable on $(0, T)$, where $u_{\le Q }$ represents the low modes of Littlewood-Paley projection of the velocity $u$.  
\end{abstract}

\maketitle

\section{Introduction}

Consider the three-dimensional incompressible Boussinesq equations, 
\begin{align} 
    \label{eq:BousU}
    &\partial_t u + (u \cdot \nabla)  u = -\nabla p  + \nu \Delta u +\theta e_3, \\
    \label{eq:BousTheta}
    & \partial_t \theta + (u \cdot \nabla)  \theta =  \kappa \Delta \theta, \\
    \label{eq:BousDiv}
    &\nabla \cdot u = 0,
\end{align}
with initial data
\begin{align}
    \label{eq:InitalU}
    &u(x,0) = u_0 (x), \\
    \label{eq:InitialTheta}
    & \theta(x, 0) = \theta_0(x),
\end{align}
where $x \in \Real^3$, $t \ge 0$, $u = u(x, t)$ is the velocity vector with divergence-free initial data, $p = p(x, t)$ is the pressure scalar, and $\theta = \theta(x, t)$ is the temperature scalar. The fluid kinematic viscosity is $\nu \ge 0$, the thermal diffusivity is $\kappa \ge 0$, and $e_3 = (0, 0, 1)^{T}.$  When $\theta$ vanishes, the system reduces to the incompressible Navier-Stokes equations, which can be further reduced to the incompressible Euler equations by setting $\nu = 0$.

The Boussinesq equations are important in the study of atmospheric sciences and they yield a wealth of interesting and difficult problems from a mathematical perspective.  The regularity of \eqref{eq:BousU}-\eqref{eq:BousDiv} has been studied throughly on its own and in relation to the regularity of other equations, such as the Navier-Stokes equations, Euler equations, and magneto-hydrodynamics (MHD) equations. 

Regularity criteria for the Boussinesq and related equations mentioned above can be split into different classes, one of which is Beale-Kato-Majda-like (abbreviated ``BKM-like") criteria. The original result by Beale, Kato, and Majda \cite{BKM} placed a condition on the vorticity in the Navier-Stokes equations. They proved if 
\begin{align}
    \label{eq:BKM}
    \int_0^T \| \nabla \times u \|_{L^\infty} ~\dt < \infty,
\end{align}
then a smooth solution to the Navier-Stokes equations on $(0, T)$ does not blow up at time $T$. This condition was weakened for the Euler equations by Planchon \cite{Planchon} and improved for the three-dimensional Navier-Stokes equations by Cheskidov and Shvydkoy \cite{CheskidovShvydkoy}.  In \cite{CheskidovDai}, Cheskidov and Dai developed BKM-like, but weaker, regularity criterion for the three-dimensional MHD equations.  

A related family of regularity criteria fall into the Ladyzhenskaya-Prodi-Serrin category. For the Navier-Stokes equations, the condition is 
\begin{align} 
	\label{eq:LPS}
	u \in L^p(0, T; L^r), ~~\mathrm{for}~~ \frac{2}{p} + \frac{3}{r} = 1, ~~r \in (3, \infty].
\end{align}
Various extensions and improvements of this type of criteria have been developed since then, such as the extension to the case for $r=3$ by Escauriaza, Seregin, and \v{S}ver\'ak  \cite{EskSerShver} and extensions to Besov spaces.

Both kinds of regularity criteria were developed for the three-dimensional incompressible Boussinesq equations.  In \cite{QiuDuYaoB} and \cite{QiuDuYaoS}, Qiu, Du, and Yao developed Serrin-type blow-up criteria for the Boussinesq equations. Particularly in \cite{QiuDuYaoB}, the authors showed a smooth solution to  \eqref{eq:BousU}-\eqref{eq:BousDiv} on time interval $[0, T)$ will remain smooth at time $T$ if $u \in L^q\big(0, T; B^s_{p, \infty}(\Real^3) \big)$ for $\frac{2}{q} + \frac{3}{p} = 1+s$, $\frac{3}{s+1}<p \le \infty$, $-1 < s \le 1$, and $(p, s) \neq (\infty, 1)$.  Ishimura and Morimoto \cite{IshimuraMorimoto} proved the Beale-Kato-Majda-like regularity criterion $\nabla u \in L^1\big(0, T; L^\infty(\Real^3)\big)$. Later Fan and Zhou \cite{FanZhou} studied the Boussinesq equations with partial viscosity and proved BKM-like regularity criteria in terms of the vorticity: $ \nabla \times u \in L^1 \big(0, T; \dot{B}^0_{\infty, \infty} (\Real^3) \big)$.  More regularity conditions in the three-dimensional case were developed in the following years, some of which can be found in \cite{QinYangWangLiu}, \cite{Xiang}, \cite{XuZhangZheng}, \cite{YangZhang}, \cite{Ye}, and \cite{Zhang}.  

A great deal of literature has also been produced for the two-dimensional case. We will not discuss these results, but rather refer the reader a few publications  on that topic: \cite{AbidiHmidi}, \cite{Chae}, \cite{ChaeKimNam}, \cite{ChaeNam}, \cite{HmidiKeranniD}, \cite{HmidiKeraaniV}, \cite{HouLi}, \cite{Taniuchi}, and \cite{Xu}.  

The main theorem discussed here also falls into the two main kinds of regularity criteria discussed above. We will show:
\begin{thm}
    \label{thm:Main}
    Let $\big(u, \theta \big)$ be a weak solution to \eqref{eq:BousU}-\eqref{eq:BousDiv} on $[0, T]$, assume $\big(u, \theta \big)$ is regular on $(0, T)$, and 
    \begin{align} 
         \| u_{\le Q} \|_{B^1_{\infty, \infty}}  \in L^1(0,T).
    \end{align} 
    Then $\big( u(t), \theta(t) \big)$ is regular on $(0,T]$.
\end{thm}

\begin{remark}
We note that the above BKM-like regularity criterion also recovers the previous known Prodi-Serrin-type regularity, in particular we improve upon the results in \cite{QiuDuYaoB}, by recovering the whole range, including the endpoint $(p, s) = (\infty, 1)$. Further, the criterion in Theorem \ref{thm:Main} improves previous BKM-like criterion for the three-dimensional Boussinesq equations since we only impose a condition on a finite amount of modes of the projection of the velocity $u$. 
\end{remark}

\begin{remark}
The notation $u_{\le Q}$ denotes the low modes of $u$ and $B^1_{\infty, \infty}$ is a Besov space. More precise definitions are presented in the next section (see \eqref{eq:LowHighNotation} and Definition \ref{def:Besov}). 
\end{remark}


\section{Background}  

\subsection{Littlewood-Paley Decomposition} 
We utilize Littlewood-Paley decomposition in our proof.  Denote wave numbers as $\la_q = 2^q$ (in some wave units). For $\psi \in C^\infty(\Real^3)$, define 
\begin{displaymath}
   \psi(\xi) = \left\{
     \begin{array}{ll}
       1, &\mathrm{for} ~~ | \xi | \le \frac{3}{4} \\
       0, &\mathrm{for} ~~ | \xi | \ge 1.
     \end{array}
   \right.
\end{displaymath} 
Next define $\phi(\xi) = \psi(\xi / \la_1) - \psi(\xi)$ and
\begin{displaymath}
   \phi_q(\xi) = \left\{
     \begin{array}{ll}
       \phi(\xi / \la_q), &\mathrm{for}~~ q \ge 0 \\
       \psi(\xi), &\mathrm{for}~~ q = -1.
     \end{array}
   \right.
\end{displaymath} 
The Littlewood-Paley projection operator $\Delta_q$ is defined as 
\begin{align}
    \Delta_q  u =  \int_{\Real^3} u(x-y) \mathcal{F}^{-1}(\phi_q)(y) \dy,
\end{align}
where $\mathcal{F}^{-1}$ is the inverse Fourier transform. Primarily, we will denote $\Delta_q u$, the  $q^{th}$ Littlewood-Paley piece of $u$, as $ u_q$ instead. 
 In the sense of distributions, one has  
\begin{align}
    \label{eq:LPDecomp}
    u = \sum_{q = -1}^\infty u_q.
\end{align}
 We also define
\begin{align}
	\label{eq:LowHighNotation}
    u_{\le Q} = \sum_{ q = -1}^Q u_q,  ~ u_{\ge Q} = \sum_{q = Q}^\infty u_q,
\end{align}
and 
\begin{align*}
    \tilde{u}_q = u_{q-1} + u_q + u_{q+1},
\end{align*}
which will be useful notation for the proof. 

\subsection{Notation, Spaces, and Solutions}
We will use the symbol $\lesssim$ (or $\gtrsim$) to mean that an inequality holds up to an absolute constant, we will denote $L^p$-norms as $\| \cdot \|_p$, and $( \cdot, \cdot)$ will denote the $L^2$ inner product.  

We use Littlewood-Paley decomposition to define some useful norms. Regularity (see Definition \ref{def:reg} ) is defined via Sobolev norms. 

\begin{definition} The homogeneous Sobolev space $\dot{H}^s$ contains functions $u$ such that the associated norm
\begin{align} 
    \label{eq:SobNormDef}
    \| u\|_{\dot{H}^s} = \Big( \sum_{q = -1}^\infty \la_q^{2s} \|u_q\|_{2}^2  \Big)^\frac{1}{2},
\end{align}
for $s \in \Real$, is finite. Note that $\| u\|_{H^s} \sim  \| u\|_{\dot{H}^s}$, which we will keep in mind throughout the work below. 
\end{definition}

The regularity criterion is defined in terms of the following Besov norm:
\begin{definition}\label{def:Besov} The norm of $u$ in Besov space $B^1_{\infty, \infty}$ is defined as
\begin{align}
    \label{eq:BesovNormDef}
     \| u (t) \|_{B^1_{\infty, \infty}}  = \sup_{q \ge -1 } \la_q \| u_q(t) \|_\infty.
\end{align}
The Besov space $B^1_{\infty, \infty} $ is the space of tempered distributions $u$ such that  $\| u (t) \|_{B^1_{\infty, \infty}}$ is finite.
\end{definition}

We work in the class of weak solutions:
\begin{definition}
A weak solution of \eqref{eq:BousU}-\eqref{eq:BousDiv} on $[0,T]$ is a pair of functions $(u,\theta)$, $u$ divergence free, in the class
\begin{align*}
    u, \theta \in C_w \big( [0,T]; L^2(\Real^3) \big) \cap L^2 \big( 0,T; H^1(\Real^3) \big)
\end{align*}
such that
\begin{multline*}
    \big( u(t), \phi(t) \big) - \big( u_0, \phi(0) \big) \\
    = \int_0^t \big(u(s), \partial_s \phi(s) \big)
    + \nu \big(u(s), \Delta \phi(s) \big) 
    + \big( u(s) \cdot \nabla \phi(s) , u(s)  \big)
    + \big ( \theta(s) e_3, \phi (s) \big) ~\ds
\end{multline*}
and
\begin{multline*}
    \big( \theta(t), \phi(t) \big) - \big( \theta_0, \phi(0) \big) \\
    = \int_0^t \big( \theta(s), \partial_s \phi(s) \big)
    + \kappa \big( \theta(s), \Delta \phi(s) \big) 
    + \big( u(s) \cdot \nabla \phi(s) , \theta(s)  \big) ~\ds,
\end{multline*}
for all divergence free test functions $\phi \in C_0^\infty \big( [0, T] \times \Real^3 \big)$. 
\end{definition}

\begin{definition}
\label{def:reg}
A Leray-Hopf weak solution of  \eqref{eq:BousU}-\eqref{eq:InitialTheta} is regular on time interval $\mathcal{I}$ if the Sobolev norm $ \| u \|_{H^s}$ is continuous for $s \ge \frac{1}{2}$ on $\mathcal{I}$.
\end{definition}

\begin{remark}
One can apply a standard bootstrap argument to show if a solution is regular, then $u$ and $\theta$ are smooth. 
\end{remark}

\subsection{The Dissipation Wave Number}
The development of our regularity criterion is linked to Kolmogorov's theory of turbulence and the dissipation wave number. The dissipation wave number is a time-dependent function that separates the low frequency inertial range, where the nonlinear term dominates the dynamics, from the high frequency dissipative range, where viscous forces takeover. In \cite{CheskidovShvydkoy}, Cheskidov and Shvydkoy defined the dissipation wave number and proved that if it belongs to $L^{5/2}(0,T)$, then the solution of the Navier-Stokes equations is regular up to time $T$. They also showed that the dissipation wave number belongs to $L^1(0,T)$ for every Leray-Hopf solution. 

We define the dissipation wave number $\Lambda(t)$ for the Boussinesq equations in an analogous manner:
\begin{align}
	& Q(t) = \min \{ q : \la_p^{-1} \|u_p \|_\infty < c \min\{\nu, \kappa \}, \forall p>q, q \ge 0 \}, \\
	& \Lambda(t) = \la_{Q(t)},
\end{align}
for absolute constant $c>0$. 

Work with the dissipation wave number and determining modes have provided key improvements to previous known regularity results for the surface quasi-geostrophic equations and the magneto-hydrodynamics equations (see \cite{CheskidovDaiSQG}, \cite{CheskidovDai}, and \cite{CheskidovDaiKavlie}), as well.

\section{Proof of the Main Theorem}
We prove Theorem \ref{thm:Main} in this section.

\begin{proof}
We test \eqref{eq:BousU} with $\Delta_q^2 u$ and \eqref{eq:BousTheta} with  $\Delta_q^2 \theta$, which yields 
\begin{align}
    \label{eq:BousUTested}
    &\frac{1}{2}\ddt \|u_q\|_2^2 \le -\nu \| \nabla u_q \|_2^2  
    + \int_{\Real^3} \Delta_q \big( u \cdot \nabla u \big) \cdot u_q \dx 
    - \int_{\Real^3} \Delta_q \big( \theta e_3 \big) \cdot u_q \dx, \\
    \label{eq:BousThetaTested}
    &\frac{1}{2}\ddt \| \theta_q\|_2^2 \le -\kappa \| \nabla \theta_q \|_2^2  
    + \int_{\Real^3} \Delta_q \big( u \cdot \nabla \theta \big) \cdot \theta_q \dx,
\end{align} 
and then multiply \eqref{eq:BousUTested} by $\la_q^{2s}$and \eqref{eq:BousThetaTested} by $\la_q^{2\sigma}$, add those two equations together, and sum over $q$ to arrive at 
\begin{multline}
    \label{eq:SobolevInequality}
    \frac{1}{2} \ddt \sum_{q = -1}^\infty \big( \la_q^{2s} \| u_q \|_2^2 + \la_q^{2\sigma}  \| \theta_q \|_2^2 \big) 
    \le - \sum_{q = -1}^\infty  \big( \la_q^{2s} \nu \| \nabla u_q \|_2^2 + \la_q^{2\sigma} \kappa \| \nabla \theta_q \|_2^2 \big) 
    \\ + I_1 + I_2 +I_3, 
\end{multline}
where
\begin{align}
    \label{eq:Int1}
    &I_1 = \sum_{q = -1}^\infty \la_q^{2s} \int_{\Real^3} \Delta_q \big( u \cdot \nabla u \big) \cdot u_q \dx, \\
    \label{eq:Int2}
    &I_2 = -\sum_{q = -1}^\infty \la_q^{2s} \int_{\Real^3} \Delta_q \big( \theta e_3 \big) \cdot u_q \dx, \\
    \label{eq:Int3}
    &I_3 = \sum_{q = -1}^\infty \la_q^{2\sigma}  \int_{\Real^3} \Delta_q \big( u \cdot \nabla \theta \big) \cdot \theta_q \dx.
\end{align} 

For \eqref{eq:Int3}, we use a similar method as in \cite{CheskidovDai}. First, we decompose \eqref{eq:Int3} into three parts:
\begin{align}
    \label{eq:Int3Pieces}
    I_3 = &\sum_{q = -1}^\infty \sum_{|q-p| \le 2} \la_q^{2\sigma} \int_{\Real^3} \Delta_q (u_{\le p-2} \cdot \nabla \theta_p) \theta_q ~\dx  \notag \\
    &+ \sum_{q = -1}^\infty \sum_{|q-p| \le 2} \la_q^{2\sigma} \int_{\Real^3}  \Delta_q ( u_p \cdot \nabla \theta_{\le p-2}) \theta_q ~\dx \\
    &+  \sum_{q = -1}^\infty \sum_{p \ge q-2}  \la_q^{2\sigma} \int_{\Real^3} \Delta_q (u_p \cdot \nabla \tilde{\theta}_{p}) \theta_q ~ \dx \notag \\
    =& I_{3,1} + I_{3,2}+ I_{3,3} \notag.
\end{align}
We use H{\"o}lder's inequality on $I_{3,2}$ to find 
\begin{align*}
    | I_{3,2} | 
    \lesssim  &\sum_{q = -1}^\infty  \la_q^{2\sigma} \| u_q \|_\infty \sum_{|q-p| \le 2} \| \theta_p \|_2 \sum_{p' \le p-2} \la_{p'} \| \theta_{p'} \|_2. \\
\end{align*}
Then we split the sum into high and low modes. For the high modes we use the definition of $\Lambda(t)$ and for the low modes we use 
\begin{align}
    \label{eq:f(t)}
    f(t) = \| u_{\le Q(t)}(t) \|_{B^1_{\infty, \infty}}  = \sup_{q \le Q(t)} \la_q \| u_q(t) \|_\infty,
\end{align}
(which will be used to define the regularity criterion later) to find
\begin{align*}
    | I_{3,2} | 
    \lesssim &c \kappa \sum_{q = Q+ 1}^\infty \la_q^{2\sigma+1} \sum_{|q-p| \le 2} \| \theta_p \|_2 \sum_{p' \le p-2} \la_{p'} \| \theta_{p'} \|_2 \\
    &+ f(t) \sum_{q= -1}^{ Q} \la_q^{2\sigma-1}  \sum_{|q-p| \le 2} \| \theta_p \|_2 \sum_{p' \le p-2} \la_{p'} \| \theta_{p'} \|_2 \\
     \lesssim &c \kappa \sum_{q = Q-1}^\infty \la_q^{2\sigma+1} \| \theta_q \|_2 \sum_{p' \le q} \la_{p'} \| \theta_{p'} \|_2 \\
    &+f(t)  \sum_{q=-1}^{Q+2} \la_q^{2\sigma-1} \| \theta_q \|_2  \sum_{p' \le q} \la_{p'} \| \theta_{p'} \|_2. \\
\end{align*}
After a rearrangement and an application of Jensen's inequality, we arrive at 
\begin{align*}
    | I_{3,2} | 
    \lesssim &c \kappa  \sum_{q = Q-1}^\infty \la_q^{\sigma+1} \| \theta_q \|_2 \sum_{p' \le q} \la_{q-p'}^{\sigma} \la_{p'}^{\sigma+1} \| \theta_{p'} \|_2 \\
    &+f(t) \sum_{q= -1}^{Q+2} \la_q^{\sigma} \| \theta_q \|_2 \sum_{p' \le q} \la_{q-p'}^{\sigma-1} \la_{p'}^{\sigma} \| \theta_{p'} \|_2 \\
    \lesssim & c \kappa \sum_{q = -1}^\infty \la_q^{2\sigma + 2} \| \theta_q \|_2^2 + f(t) \sum_{q = -1}^{Q+2} \la_q^{2\sigma} \| \theta_q \|_2^2,
\end{align*}
for 
\begin{align}
	\label{eq:sigmaUpperbound}
	\sigma < 0. 
\end{align}

By Bony's paraproduct and commutator notation, which says 
\begin{align*}
    [ \Delta_q, u_{\le p-2} \cdot \nabla ] \theta_p  = \Delta_q (u_{\le p -2} \cdot \nabla \theta_p) - u_{\le p -2} \cdot \nabla \Delta_q \theta_p,
\end{align*}
one may decompose $I_{3,1}$ as
\begin{align}
    \label{eq:I31Pieces}
    I_{3,1}  
    =& \sum_{q = -1}^\infty \sum_{|q-p| \le 2} \la_q^{2 \sigma} \int_{\Real^3} [\Delta_q, u_{\le p-2} \cdot \nabla] \theta_p \theta_q ~\dx  \notag \\
    & + \sum_{q = -1}^\infty \la_q^{2\sigma} \int_{\Real^3} u_{\le q-2} \cdot  \nabla \theta_q \theta_q ~\dx \\
    & + \sum_{q = -1}^\infty \sum_{|q-p| \le 2} \la_q^{2\sigma} \int_{\Real^3} (u_{\le p-2} - u_{\le q-2}) \cdot \nabla \Delta_q \theta_p \theta_q   ~\dx  \notag \\
    =& I_{3,1,1} + I_{3,1,2}+I_{3, 1, 3} .    \notag
\end{align}
In \cite{CheskidovDai}, they note that their term (equivalent to our $I_{3,1,1}$) can be estimated as
\begin{align*}
    | I_{3,1,1} | \lesssim c \kappa \sum_{q = Q+2}^\infty \la_q^{2\sigma+2} \|\theta_q \|_2^2 
    + f(t) \sum_{q = -1}^\infty \la_q^{2\sigma} \| \theta_q \|_2^2.
\end{align*}
The second term, $I_{3,1,2} =  0$ because of the divergence-free condition on $u$. We also refer the reader to \cite{CheskidovDai}, where one can find
\begin{align*}
    | I_{3, 1, 3}| + |I_{3,3}| \lesssim c \kappa \sum_{q = -1}^\infty \la_q^{2\sigma+2} \|\theta_q \|_2^2 
    + f(t) \sum_{ q = -1}^{Q+2} \la_q^{2\sigma} \| \theta_q \|_2^2.
\end{align*}  

The above estimates on the pieces of \eqref{eq:Int3} yield 
\begin{align}
    \label{eq:EstimateInt3}
    |I_3| \lesssim c \kappa \sum_{q = -1}^\infty \la_q^{2\sigma+2} \| \theta_q\|_2^2 + f(t) \sum_{q = -1}^\infty  \la_q^{2\sigma} \| \theta_q \|_2^2.
\end{align}

For \eqref{eq:Int1}, we refer the reader to the estimates carried out in \cite{CheskidovShvydkoy} on the Navier-Stokes equations (and the equivalent term in  \cite{CheskidovDai} on the MHD equations), where they show 
\begin{align} 
    \label{eq:EstimateInt1}
    | I_1 | \lesssim c \nu \sum_{q = -1}^\infty \la_q^{2s+2} \| u_q \|_2^2 
    + f(t) \sum_{q = -1}^\infty \la_q^{2s} \|u_q \|_2^2, 
\end{align}
where $f(t)$ is the same as in \eqref{eq:f(t)}.
The bound \eqref{eq:EstimateInt1} holds for 
\begin{align}
	\label{eq:sBounds}
	\frac{1}{2} \le s < 1.
\end{align}

We use Young's inequality to estimate \eqref{eq:Int2} as
\begin{align}
    \label{eq:EstimateInt2}
    | I_2 | = \Big|  \sum_{q = -1}^\infty  \la_q^{2s} \int_{\Real^3} \Delta_q \big( \theta e_3 \big) \cdot u_q \dx \Big|  
    \lesssim \sum_{q = -1}^\infty  \la_q^{2s}  \big( \| u_q \|_2^2 + \| \theta_q \|_2^2 \big).
\end{align}
We may break up this sum as follows:
\begin{align}
    \label{eq:EstimateInt2Broken}
    | I_2 | \lesssim \sum_{q = -1}^\infty  \la_q^{2s} \| u_q \|_2^2 +  \sum_{q= -1}^N \la_q^{2s}  \| \theta_q \|_2^2 + \sum_{q =N }^\infty \la_q^{2s}  \| \theta_q \|_2^2,
\end{align}
where $N$ is finite and will be determined later.

Estimates \eqref{eq:EstimateInt1}, \eqref{eq:EstimateInt2Broken}, and \eqref{eq:EstimateInt3} in \eqref{eq:SobolevInequality} yield
\begin{align}
    \label{eq:SobolevIneq2}
    \frac{1}{2} \ddt \sum_{q = -1}^\infty  \big( \la_q^{2s} \|u_q\|_2^2 + \la_q^{2\sigma} \| \theta_q \|_2^2 \big) 
     \lesssim & - \nu \sum_{q = -1}^\infty  \la_q^{2s+2 } \| u_q \|_2^2 + c \nu \sum_{q = -1}^\infty \la_q^{2s+2} \| u_q \|_2^2 \notag \\
     & - \kappa \sum_{q = -1}^\infty \la_q^{2\sigma +2}  \| \theta_q \|_2^2  + c \kappa \sum_{q = -1}^\infty \la_q^{2\sigma +2} \| \theta_q \|_2^2  \notag \\
     & +  \sum_{q= -1}^N \la_q^{2s}  \| \theta_q \|_2^2 + \sum_{q =N }^\infty \la_q^{2s}  \| \theta_q \|_2^2   \\
     & + \big( f(t) + 1 \big) \sum_{q = -1}^\infty \la_q^{2s} \|u_q\|_2^2  \notag \\
     & + f(t) \sum_{q = -1}^\infty \la_q^{2\sigma} \| \theta_q \|_2^2. \notag
\end{align}
For 
\begin{align}
	\label{eq:sANDsigma}
	2s < 2\sigma +2
\end{align}	
and absolute constants $C_1, C_2, C_3, C_4$ and $C_5$, we have 
\begin{align}
    \label{eq:SobolevIneq3}
    \frac{1}{2} \ddt \sum_{q = -1}^\infty  \big( \la_q^{2s} \|u_q\|_2^2 + \la_q^{2\sigma} \| \theta_q \|_2^2 \big) 
     \le & - \nu \sum_{q = -1}^\infty  \la_q^{2s+2 } \| u_q \|_2^2 + C_1 c \nu \sum_{q = -1}^\infty \la_q^{2s+2} \| u_q \|_2^2 \notag \\
     & - \kappa \sum_{q = -1}^\infty \la_q^{2\sigma +2}  \| \theta_q \|_2^2  + C_2 c \kappa \sum_{q = -1}^\infty \la_q^{2\sigma +2} \| \theta_q \|_2^2  \notag \\
     & +  C_3 \sum_{q= -1}^N \la_q^{2s}  \| \theta_q \|_2^2 + C_4 \la_N^{2s-2\sigma-2}\sum_{q =N }^\infty \la_q^{2\sigma+2}  \| \theta_q \|_2^2   \\
     & + C_5  \big( f(t) + 1 \big) \sum_{q = -1}^\infty \big( \la_q^{2s} \|u_q\|_2^2  + \la_q^{2\sigma} \| \theta_q \|_2^2 \big).  \notag
\end{align}
For $N$ large enough, one may choose $N = N(\kappa)$ such that $C_4 \la_N^{2s -2\sigma -2} \le \kappa/2$. In fact, one may solve for this $N$ explicitly. Since such a finite $N$ exists, then we can use the following two facts:
\begin{align}
	\label{eq:BelowN}
	\sum_{q= -1}^N \la_q^{2s}  \| \theta_q \|_2^2 < \infty,
\end{align}
and 
\begin{align}
	\label{eq:AboveN}
	C_4 \la_N^{2s-2\sigma-2}\sum_{q =N }^\infty \la_q^{2\sigma+2}  \| \theta_q \|_2^2 \le \frac{\kappa}{2} \sum_{q = -1}^\infty \la_q^{2\sigma+2}  \| \theta_q \|_2^2.
\end{align}
Using \eqref{eq:BelowN} and \eqref{eq:AboveN} in \eqref{eq:SobolevIneq3}, we arrive at the following differential inequality:
\begin{align}
    \label{eq:SobolevIneq4}
    \frac{1}{2} \ddt \sum_{q = -1}^\infty  \big( \la_q^{2s} \|u_q\|_2^2 + \la_q^{2\sigma} \| \theta_q \|_2^2 \big) 
     \le & (- \nu +C_1 c \nu) \sum_{q = -1}^\infty  \la_q^{2s+2 } \| u_q \|_2^2  \notag \\
     & (- \frac{\kappa}{2} +C_2 c \kappa) \sum_{q = -1}^\infty \la_q^{2\sigma +2}  \| \theta_q \|_2^2  \notag \\
     & + C_5  \big( f(t) + 1 \big) \sum_{q = -1}^\infty \big( \la_q^{2s} \|u_q\|_2^2  + \la_q^{2\sigma} \| \theta_q \|_2^2 \big) \\
     & + C_3 \sum_{q= -1}^N \la_q^{2s}  \| \theta_q \|_2^2 \notag 
\end{align}

The choice $c = \min\{\frac{1}{C_1}, \frac{1}{2C_2} \}$ yields  
\begin{multline}
    \label{eq:SobolevIneq4}
    \ddt \big( \|u\|_{\dot{H}^s}^2 +  \| \theta \|_{\dot{H}^\sigma}^2  \big) 
    \le   C(\nu, \kappa, s, \sigma) \big( f(t) + 1\big) \big( \|u\|_{\dot{H}^s}^2 +  \| \theta \|_{\dot{H}^\sigma}^2  \big) \\ 
     + C_6\sum_{q= -1}^N \la_q^{2s}  \| \theta_q \|_2^2.
\end{multline}
By Gr\"onwall's inequality (noting \eqref{eq:BelowN}), we have that $\|u\|_{\dot{H}^s}^2 +  \| \theta \|_{\dot{H}^\sigma}^2 $ remains bounded on $(0, T)$ for $\frac{1}{2} \le s < 1$ and $s-1 < \sigma < 0$ if 
\begin{align}
    \label{eq:Criterion}
    \| u_{\le Q}\|_{B^1_{\infty, \infty}}   \in L^1(0, T).
\end{align}
Thus, by Definition \ref{def:reg}, we reach the desired conclusion. 
\end{proof}

\textbf{Acknowledgments:}  The work of K. Zaya was supported by NSF grants DMS-1210896 and DMS-1517583 while at the University of Illinois at Chicago, and later by NSF grant DMS-1515161 while at the University of Michigan.

\bibliographystyle{abbrv}
\bibliography{BoussinesqRegularityCriteria}

\end{document}